\documentclass[12pt]{article}
\usepackage{amsmath}
\usepackage{amsfonts}
\usepackage{amssymb}
\usepackage{exscale,relsize}
\usepackage{amsthm}
\usepackage{hyperref}
\usepackage{graphicx}
\usepackage{color}
\usepackage{subfig}
\usepackage[mathscr]{eucal}

\def\li{\rm {li}}

\def\DHLhksqrt#1#2{\setbox0=\hbox{$#1\sqrt{#2\,}$}\dimen0=\ht0
\advance\dimen0-0.2\ht0
\setbox2=\hbox{\vrule height\ht0 depth -\dimen0}%
{\box0\lower0.4pt\box2}}
\def\bee{\begin{equation}}
\def\eee{\end{equation}}

\def\DHLhksqrt#1#2{\setbox0=\hbox{$#1\sqrt{#2\,}$}\dimen0=\ht0
\advance\dimen0-0.2\ht0
\setbox2=\hbox{\vrule height\ht0 depth -\dimen0}%
{\box0\lower0.4pt\box2}}

\allowdisplaybreaks

\begin{document}

\thispagestyle{empty}
\bigskip
\centerline{    }
\centerline{\Large\bf  6+$\infty$ new  expressions}
\centerline{             }
\centerline{\Large\bf for the  Euler-–Mascheroni constant}
\bigskip\bigskip
\centerline{\large\sl Marek Wolf}
\begin{center}
Cardinal  Stefan  Wyszynski  University, Faculty  of Mathematics and Natural Sciences. College of Sciences,\\
ul. W{\'o}ycickiego 1/3,   PL-01-938 Warsaw,   Poland, e-mail:  m.wolf@uksw.edu.pl
\end{center}

\begin{abstract}
In the first part we present results of  four  ``experimental''  determinations of the  Euler–-Mascheroni constant   $\gamma$.
Next we give new formulas  expressing the $\gamma$ constant  in terms of the  Ramanujan--Soldner constant  $\mu$. Employing  the
cosine  integral  we  obtain  the infinity  of formulas for $\gamma$.
\end{abstract}

\bigskip\bigskip\bigskip

\bibliographystyle{abbrv}

\section{Introduction}
\label{wstep}
The Euler--Mascheroni constant is defined by the following limit:
\bee
\gamma = \lim_{k \rightarrow \infty}  \left( \sum_{n=1}^k \frac{1}{n} - \log(k) \right)=0.57721566490153286\ldots.
\label{gamma}
\eee
see  \cite{Havil-2003},  \cite{Lagarias_2013}.  It is  not  known whether $\gamma$  is irrational, see \cite{Sondow_2003},
\cite{Lagarias_2013}.     The limit in \eqref{gamma} is very slowly convergent (like $n^{-1}$)  and   in
\cite{detemple1993} it was  shown that  slight modification of  \eqref{gamma}:
\[
\gamma=\lim_{k \to \infty} \left( \sum_{n=1}^k \frac{1}{n} - \log(k+\frac{1}{2}) \right)
\]
improves convergence  to $1/n^2$.  Presently  sequences  converging to $\gamma$  much  faster  are known, see
\cite{Mortici_2010}  where sequence which converge to $\gamma$  like $n^{-6}$  is presented.   There is a lot
of formulas expressing $\gamma$  as  series,  integrals or products,
see  \cite{Havil-2003} and  e.g.  \cite{Gourdon_Sebah}, \cite{Choi_Srivastava_2010}, \cite{Lagarias_2013}.
In particular there  is  infinity of formulas for $\gamma$,   we mention here   \cite[p. 4]{Gourdon_Sebah}:
\bee
\gamma=\sum_{k=1}^n \frac{1}{k}- \log(n)-\sum_{k=2}^\infty \frac{\zeta(k, n+1)}{k}, ~~~~n=2, 3, \ldots~,
\eee
where the Hurwitz  zeta function:
\bee
\zeta(s, k)=\sum_{n=0}^\infty \frac{1}{(n+k)^s},   ~~~~~~~\Re(s)>1.
\eee
Another infinite  set of formulas for $\gamma$ we found in \cite[eq.(9.3.10)]{Boros_Moll}:
\bee
\gamma=\sum_{k=1}^n \frac{1}{k}-\log n -\int_n^\infty \frac{\{x\}}{x^2} dx, ~~~~n=1, 2, 3,  \ldots~.
\eee
Here $\{x\}$   is the fractional part of $x$.

There are even uncountably many formulas for $\gamma$,  see  e.g. \cite{Brent_McMillan_1980}:  for real $r>0$
\bee
\gamma=\lim_{n \to \infty} \frac{\sum_{k=0}^\infty (\frac{n^k}{k!})^r(\sum_{j=1}^k\frac{1}{j}-\log(k))}{\sum_{k=0}^\infty (\frac{n^k}{k!})^r}.
\eee
There are also  doubly uncountably formulas for $\gamma$,  we present the formula (3.13)  from \cite{Choi_Srivastava_2010}:
\bee
\gamma=r\int_0^\infty \left(\frac{1}{1+x^q}-\exp(-x^r)\right) \frac{dx}{x}, ~~~~~~~ q>0, ~r>0.
\eee
The  numerical   value of the  Euler--Mascheroni constant  was  calculated in the past many times see e.g. \cite{Brent_McMillan_1980},
the  present day world  record  is 477,511,832,674 	decimal digits of $\gamma$  and belongs to Ron Watkins, see
\verb+http://www.numberworld.org/digits/EulerGamma+

The  Euler--Mascheroni constant  appears  also in many places in number theory  and
in the  theory of the Riemann zeta  function,  for example in the Nicolas' and Robin's  criterions  for the Riemann Hypothesis,
see  e.g. \cite[vol.1, chapters 5 and 7]{Broughan_2017}.
One of the most amazing  appearances of the $\gamma$  constant
is in the F. Merten's  two  products over   primes  \cite[p.351]{H-W},,  one of which involves  constants
$\pi, ~e, ~\gamma$  (``holy   trinity''):
\bee
\lim_{n\to \infty}  \frac{1}{ \log(n)}\prod_{p<n} \Big(1+\frac{1}{p}\Big) = \frac{6e^\gamma}{\pi^2}
\label{Mertens}
\eee
from which we obtain
\bee
\gamma=\log\Big(\frac{\pi^2}{6}\lim_{n\to \infty}  \frac{1}{ \log(n)}\prod_{p<n} \Big(1+\frac{1}{p}\Big)\Big).
\eee
With present day computers we can check the accuracy of the  above relation.  In  the Table I we present lhs and  rhs  of
\eqref{Mertens}  as well as  computed from the values of the  finite  products over primes values  of the
\bee
\gamma(n)=\log\Big(\frac{\pi^2}{6} \frac{1}{ \log(n)}\prod_{p<n} \Big(1+\frac{1}{p}\Big)\Big).
\eee

\vskip 0.4cm
\begin{center}
{\sf TABLE 1} ~The values of the  product in \eqref{Mertens} up to $n=1000, 10000, \ldots, 10^{13}$  (second column)
and values of its values following from the Mertens formula (third column),  theirs ratio in fourth column and  finite approximations
to $\gamma$ in the last column.  The fluctuations in the last digits of the values obtained from the computer
are presumably caused by   cumulation the floating--point errors.\\
\bigskip
\begin{tabular}{|c|c|c|c|c|} \hline
$n$ & $ \prod_{p<n} (1+1/p_n) $ & $6e^{\gamma} \log(n)/\pi^2  $  &  ratio & $\gamma(n)$ \\ \hline
 $  10^{  3}  $ &        7.5094464  &           7.4891425 &     1.0027111 &    0.57992110 \\ \hline
 $  10^{  4}  $ &        9.9849904  &           9.9733461 &     1.0011675 &    0.57838053 \\ \hline
 $  10^{  5}  $ &       12.4756558  &          12.4721158 &     1.0002838 &    0.57749746 \\ \hline
 $  10^{  6}  $ &       14.9651229  &          14.9643917 &     1.0000489 &    0.57726252 \\ \hline
 $  10^{  7}  $ &       17.4570890  &          17.4568441 &     1.0000140 &    0.57722769 \\ \hline
 $  10^{  8}  $ &       19.9494269  &          19.9493052 &     1.0000061 &    0.57721977 \\ \hline
 $  10^{  9}  $ &       22.4418428  &          22.4417674 &     1.0000034 &    0.57721703 \\ \hline
 $  10^{ 10}  $ &       24.9342956  &          24.9342295 &     1.0000027 &    0.57721631 \\ \hline
 $  10^{ 11}  $ &       27.4267504  &          27.4266917 &     1.0000021 &    0.57721581 \\ \hline
 $  10^{ 12}  $ &       29.9192150  &          29.9191539 &     1.0000020 &    0.57721571 \\ \hline
 $  10^{ 13}  $ &       32.4116846  &          32.4116161 &     1.0000021 &    0.57721578 \\ \hline
\end{tabular}
\end{center}

\bigskip
\bigskip

The average value of the divisor function $d(n)$ counting  the number of divisors of  $n$  including  1 and $n$ is given
by the theorem  proved by Dirichlet, see e.g. \cite[Th.320]{H-W}:
\bee
 \frac{1}{n} \sum_{k=1}^n d(k)=\log n + 2\gamma-1+\mathcal{O}\Big(\frac{1}{\sqrt{n}}\Big).
\label{divisors}
\eee
In Table 2 values of $\gamma$  obtained  from above formula for $n=2^{15}, \ldots, 2^{23}$  are  presented.

\vskip 0.4cm
\begin{center}
{\sf TABLE 2} ~The values of $\gamma$  obtained from \eqref{divisors} for $n=2^{15},  2^{16}, \ldots, 2^{23}$.\\
\bigskip
\begin{tabular}{|c|c|c|} \hline
$n$  &  $ \sum_{k=1}^n d(k)  $  & $\gamma$  from \eqref{divisors} \\ \hline
          32769  &               345785  &          0.5776565 \\ \hline
          65537  &               736974  &          0.5774880 \\ \hline
         131073  &              1564762  &          0.5773423 \\ \hline
         262145  &              3311206  &          0.5772996 \\ \hline
         524289  &              6985780  &          0.5772608 \\ \hline
        1048577  &             14698342  &          0.5772438 \\ \hline
        2097153  &             30850276  &          0.5772336 \\ \hline
        4194305  &             64607782  &          0.5772288 \\ \hline
        8388609  &            135030018  &          0.5772237 \\ \hline

 \end{tabular}
 \end{center}

 The sum  of reciprocals of non--trivial  zeros  $\rho$  of the
Riemann's  zeta  function  $\zeta(s)$  also  involves $\gamma$ \cite[p.67 and p. 159]{Edwards}:
\bee
\sum_{\rho}\frac{1}{\rho} ~=~ 2+\gamma-\log(4\pi) =  0.046191417932\ldots~.
\label{zera_rho}
\eee
The above sum is real and convergent when zeros  $\rho$ and complex conjugate $\overline{\rho}$ are  paired together
and  summed according to   increasing absolute  values of the imaginary parts  of $\rho$.    Several  years ago using the
L-function  calculator  written   by Michael Rubinstein (see  http://doc.sagemath.org/html/en/reference/lfunctions/sage/lfunctions/lcalc.html)
we have calculated  100,000,000  zeros of   $\zeta(s)$;  the last obtained zero has the value
$\rho_{100000000}=\frac12+\imath 42653549.7609515$.  In Table 2
we present approximations to $\gamma$  obtained from \eqref{zera_rho} after summing over   1000, 10,000, ..., 100,000,000  zeros
of zeta function.

The largest known  prime  numbers are of the form $\mathcal{M}_p=2^p-1$ where in turn $p$ is also a prime and they are called Mersenne primes,
see eg. https://www.mersenne.org/. In \cite[p.101]{Pomerance_1981} (see also \cite{Wagstaff1983})   the  Lenstra--Pomerance--Wagstaff
conjecture  was  formulated:   
the number of $p<x$  with $2^p-1$ prime grows  like  
\bee
\sharp \{p<x {\rm  ~~ and ~~} 2^p-1 {\rm ~~  prime  }\} \sim \frac{e^{\gamma}}{\log 2} \log x.
\label{Wagstaff}
\eee
The  presence of $\gamma$  here comes  from the Merten's  result \eqref{Mertens}.  In the Fig. 1 we compare the
Pomerance -- Wagstaff conjecture with all 51 presently known  Mersenne primes.  From the fit of actual number of $p<x$
with $2^p-1$  prime to the $\log x$  gives  rather  poor value    $\gamma=0.61$,  thus  it is rather not convenient
way to compute the  Euler--Mascheroni constant.

\vskip 0.4cm
\begin{center}
{\sf TABLE 3} The value of $\gamma$  obtained from \eqref{zera_rho} after  summing over $n=1000, 10000, \dots,  100000000$
zeros of  $\zeta(s)$.\\
\bigskip
\begin{tabular}{|c|c|} \hline
$n$ &  $\gamma$ \\ \hline
1000  &          0.5757765 \\ \hline
10000  &          0.5769463 \\ \hline
100000  &          0.5771715 \\ \hline
1000000  &          0.5772091 \\ \hline
10000000  &          0.5772147 \\ \hline
100000000  &          0.5772155 \\ \hline
\end{tabular}
\end{center}

In this paper we will present  some new formulas  for $\gamma$
obtained by  putting in the series for the logarithmic   and cosine integrals  special values for the argument.
Similar  idea  appeared in \cite{Sweeney_1963},  where  the series for the  exponential integral was  used to
calculate  $\gamma$  up to  3566 decimal  places.  A few
of these new  expressions present  the  Euler-–Mascheroni constant in the form of the difference of  two  numbers one of  which
is transcendental.   It  gives  hopes  for the proof not only of the irrationality of $\gamma$  but also its
transcendentality.

\begin{figure}
\vspace{-2.3cm}
\hspace{-3.5cm}
\begin{center}
\includegraphics[width=0.4\textheight,angle=0]{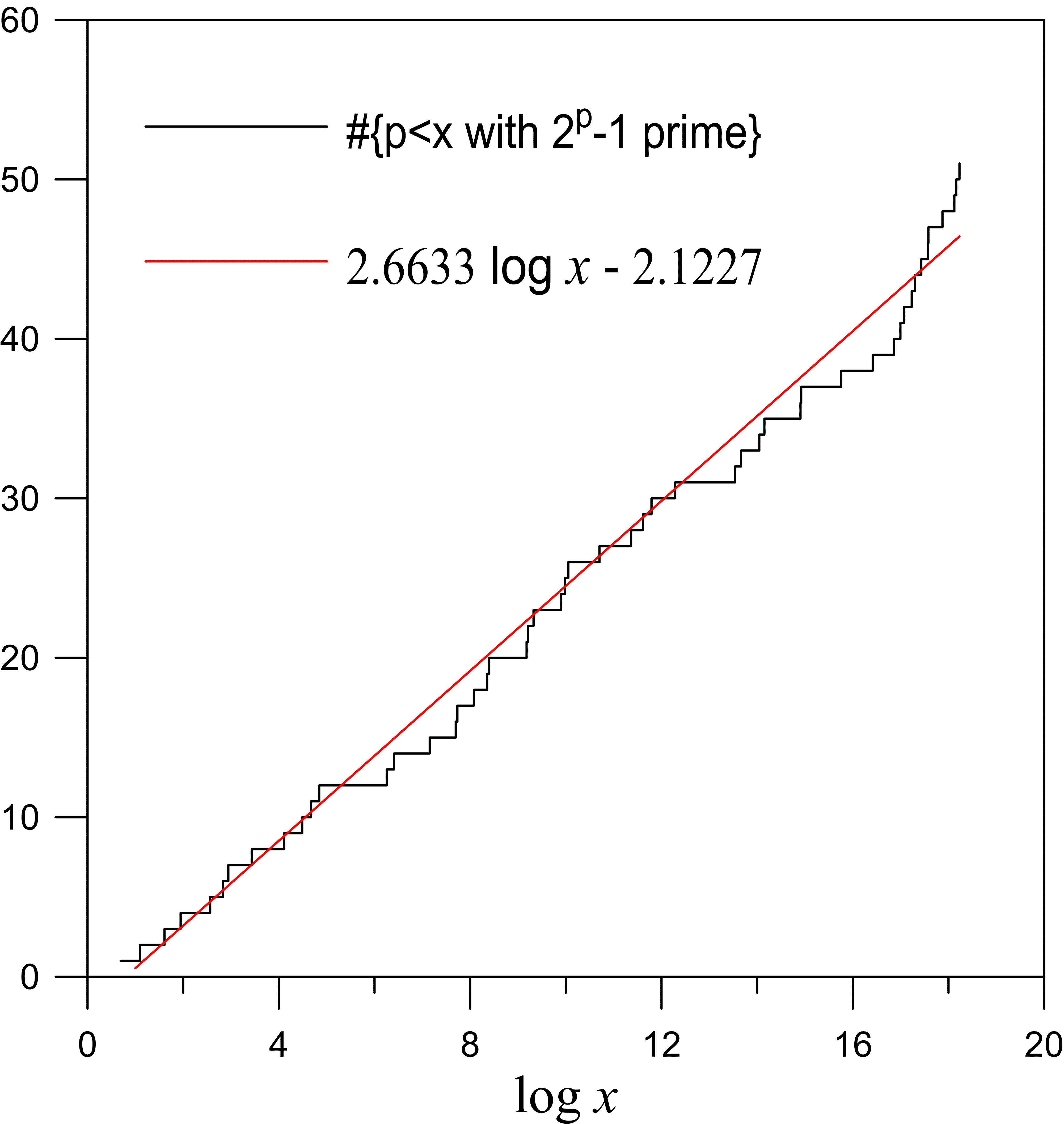} \\
\vspace{0.2cm} {\small Fig.1  The plot  illustrating  the Lenstra--Pomerance--Wagstaff conjecture.
The  least--square   fit was done  to all known $\mathcal{M}_p$ and it is
$2.6633\log x - 2.1227$, which  gives rather bad  value for $\gamma$ of 0.61}.
\end{center}
\end{figure}

\section{Logarithmic integral}
\label{LogI}

The logarithmic integral is defined for all positive real numbers  $x\neq 1$ by the definite integral
\bee
{\rm li}(x) \equiv
\left\{
\begin{array}{ll}
p.v.\displaystyle{\int_0^x \frac{du}{\log(u)},~~{\rm for}~~   x>1;}\\
\displaystyle{ \int_0^x \frac{du}{\log(u)}, ~~{\rm for}~~   0<x<1, }
\end{array}
\right.
\label{PNT}
\eee
where $p.v.$  stands for  Cauchy principal value  around $u=1$:
\bee
  p.v. \int_0^x \frac{du}{\log(u)}=\lim_{\epsilon\to 0} \Big( \int_0^{1-\epsilon}  \frac{du}{\log(u)}+ \int_{1+\epsilon}^x  \frac{du}{\log(u)} \Big).
\eee

There is the series giving logarithmic integral
${\li}(x)$ for  all $x>1$ and  quickly convergent because it  has $n!$ in denominator  and $\log^n(x)$  in  numerator
(see \cite[\S 5.1]{abramowitz+stegun})
\begin{equation}
{\li }(x) ~=~\gamma + \log \log x  + \sum_{n=1}^{\infty} {\log^{n} x \over n \cdot n!}  \quad {\rm for} ~ x > 1.
\label{Li-series}
\end{equation}
The variant of the above series after some change of variable:
\bee
\int_{x}^\infty \frac{e^{-t}}{t} dt~=~-\gamma - \log x  + \sum_{n=1}^{\infty} \frac{(-1)^{n-1} x^n}{ n \cdot n!}
\eee
was  used in \cite{Sweeney_1963} for large $x>0$, when lhs of above equation is practically  zero (in fact it is
$\mathcal{O}(e^{-x}/x)$),   to compute  3566  digits of $\gamma$,  see also  \cite{Brent_McMillan_1980}.

The  logarithmic integral takes a  value 0 at only one  number which is denoted  by  $\mu$  and is
called the Ramanujan--Soldner constant
\bee
\int_0^\mu \frac{du}{\log u} = 0,
\eee
see e.g.  \cite[eq.(11.3)]{BerndtIV} and its  numerical  value is:
\[
\mu=1.45136923488338105028396848589202745\ldots ~.
\]
Thus  for $x>\mu$  we have:
\bee
{\rm li}(x) = \int_\mu^x \frac{du}{\log(u)}.
\eee
Inserting in  \eqref{Li-series}  $x=\mu>1$ we obtain  the first  formula  expressing the Euler-–Mascheroni constant  by the
Ramanujan--Soldner  constants:     
\bee
\gamma ~=~ -\log \log \mu - \sum_{n=1}^{\infty} \frac{{\log^n} \mu}{n \cdot n!}.
\label{EMRS1}
\eee
Using PARI \cite{PARI} we  checked that summing above to  $n=20$  reproduces  31 digits of $\gamma$.   In the {\bf Appendix}
we  give the  script  to  reproduce this result  with whatever  number of digits.

Even faster converging series was discovered by Ramanujan \cite[p.130]{BerndtIV}:
\begin{equation}
\int_\mu^x \frac{du}{\log u} = \gamma + \log \log x + \sqrt{x} \sum_{n=1}^{\infty} \frac{ (-1)^{n-1} (\log  x)^n} {n! \, 2^{n-1}}
\sum_{k=0}^{\lfloor (n-1)/2 \rfloor} \frac{1}{2k+1} \quad \quad {\rm for} ~ x > 1.
\label{Li_R}
\end{equation}
Putting here $x=\mu$ we obtain  second formula  for  the  Euler–-Mascheroni constant:
\bee
\gamma ~=~ - \log\log \mu + \sqrt{\mu} \sum_{n=1}^{\infty} \frac{ (-1)^{n} (\log \mu)^n} {n! \, 2^{n-1}}
\sum_{k=0}^{\lfloor (n-1)/2 \rfloor} \frac{1}{2k+1}.
\label{EMRS2}
\eee
We checked using PARI  that summing above to  $n=20$  reproduces  correctly 37 digits of $\gamma$.

Putting in \eqref{Li-series}   $x=e$  simplifies  series and we obtain  third   expression
for the   Euler–-Mascheroni constant:
\bee
\gamma ~=~ \int_\mu^e \frac{du}{\log u} - \sum_{n=1}^{\infty} \frac{1}{n \cdot n!}:=\alpha-\beta,
\label{EM3}
\eee
where the numbers
\bee
\alpha:=\int_\mu^e \frac{du}{\log u}=1.89511781635593675546652\ldots~,
\eee
\bee
\beta:=\sum_{n=1}^{\infty} \frac{1}{n \cdot n! } = 1.31790215145440389486\ldots ~.
\label{moja_liczba}
\eee
The  number $\beta$ is  irrational as the reasoning proving the irrationality of $e=\sum_{n=0}^\infty 1/n!$
(see e.g. \cite[p.65]{Rudin-principles}) can be repeated  here {\it mutatis mutandis}. In fact from the Siegel--Shidlovsky theorem
\cite[see eq.5.2  for $k=1$]{Feldman_Shidlovskii_1967}  it follows  that \eqref{moja_liczba} is transcendental.

Putting in \eqref{Li_R}  $x=e$ we  get the  fourth  expression  for Euler--Mascheroni constant
\bee
\gamma ~=~ \int_\mu^e  \frac{du}{\log u}  + \sqrt{e} \sum_{n=1}^{\infty} \frac{ (-1)^{n}} {n! \, 2^{n-1}}
\sum_{k=0}^{\lfloor (n-1)/2 \rfloor} \frac{1}{2k+1}.
\label{EM4}
\eee

Finally let us notice  that in \cite{Havil-2003}  at several places (e.g. pp. 52, 104) we can read  that Euler had hoped that
$\gamma$   is  the logarithm of some important number. Above we have given for $\gamma$  two  series  in logarithm
of the Ramanujan--Soldner  constant.

\bigskip

\section{Cosine integral}
\label{CosI}

Many special  functions involve in theirs  expansions the  Euler-–Mascheroni constant.
The function cosine integral ${\rm Ci}(x)$  has the series expansion also containing  $\gamma$ (see  e.g.
\cite[\S5.2]{abramowitz+stegun}):
\bee
{\rm Ci}(x)~=~-\int_x^\infty \frac{\cos u}{u} du ~=~\gamma+ \log x+\sum_{n=1}^\infty \frac{(-x^2)^n}{2n(2n)!}
\label{kosinus}
\eee
\[
~~~~~~~~~~~~~~~~~~~~~~~~~~~~~~~~~~~~~~~=~\gamma+ \log x+\sum_{n=1}^\infty \frac{(-x^2)^n}{2^{n+1}n n!(2n-1)!!}
\]
because  $(2n)!=2^n n! (2n-1)!!$,  where odd  factorial $(2n-1)!!=1\cdot 3 \cdot 5 \cdot \ldots \cdot (2n-1)$.
Putting above $x=1$ we obtain the fifth expression for the  Euler-–Mascheroni constant:
\bee
\gamma~=~-\int_1^\infty \frac{\cos u}{u} du +\sum_{n=1}^\infty \frac{(-1)^{n-1}}{2n(2n)!},
\label{fifth}
\eee
where:
\bee
\int_1^\infty \frac{\cos u}{u} du = -0.3374039229009681346626\ldots
\eee
and
\bee
\sum_{n=1}^\infty \frac{(-1)^{n-1}}{2n(2n)!}~=~0.2398117420005647259439\ldots~.
\eee

The cosine integral ${\rm Ci}(x)$  has infinity of zeros  that  do not have their own names   and  are
non--periodic;  they are  usually  denoted by $c_k$, see \cite{MacLeod2002}.  The first zeros are
$c_0=0.61650548562,~c_1=3.38418042255,~ c_2=6.42704774405, ~\ldots$. A.J. MacLeod in  \cite{MacLeod2002}  gives the asymptotic
expansion  for these zeros:
\bee
c_k\sim k\pi + \frac{1}{k\pi}-\frac{16}{3(k\pi)^3}+\frac{1673}{15(k\pi)^5}-\frac{507746}{105(k\pi)^5}+\ldots.
\label{MacLeod}
\eee

Putting   zeros  $c_k$  into  \eqref{kosinus}  we obtain an  infinity of  expressions for $\gamma$
\bee
\gamma~=~-\sum_{n=1}^\infty \frac{(-c_k^2)^{n}}{2n(2n)!}-\log c_k, ~~~~~~~~k=0,1,2,\ldots~.
\label{sixth}
\eee
In the Table 4 we  present values for $\gamma$  obtained from above formula when  $c_k$  are calculated  from
\eqref{MacLeod}.  In the last  column the  difference between  values  in third column and  $\gamma$  are presented. 
\bigskip

From the MacLeod  formula  \eqref{MacLeod}  we see that  large zeros of ${\rm Ci}(x)$  approach  just zeros  of $\sin(x)=\int\cos(x)dx$:
$c_k\sim  k\pi$   for large $k$.  Thus we  have the  sixth    formula for  the   Euler-–Mascheroni constant:
\bee
\gamma = \lim_{k \to \infty} \Big(\sum_{n=1}^\infty \frac{(-1)^{n-1}(k^2 \pi^2)^{n}}{2n(2n)!}-\log(k\pi)  \Big).
\label{six}
\eee
This  formula in some sense  resembles the  original  definition \eqref{gamma}. We  checked that for $k=8000$  the expression 
in above big parentheses  gives  0.57721566648467633997998 i.e.  it reproduces  correctly  first 8 digits of $\gamma$.  
                                 
\vskip 0.4cm
\begin{center}
{\sf TABLE 4} The values of expression  \eqref{sixth} when for $c_k$ the series \eqref{MacLeod} are
substituted  for $k=10,  20 ,\ldots, 100$. \\
\bigskip
\begin{tabular}{|c|c|c|c|} \hline
$k$ & $ c_k $ from eq.\eqref{MacLeod} & eq.\eqref{sixth} for this $c_k$ & $\mid$eq.\eqref{sixth} for this $c_k-\gamma\!\!\mid$ \\ \hline
10  &  31.447589011629313  &  0.5772156649004098 & $ 1.123\times 10^{-12} $  \\ \hline
20  &  62.847747177749027  &  0.5772156649015328 & $ 1.953\times 10^{-17} $  \\ \hline
30  &  94.258383581485718  &  0.5772156649015328 & $ 2.888\times 10^{-20} $  \\ \hline
40  &  125.67166120666795  &  0.5772156649015328 & $ 2.657\times 10^{-22} $  \\ \hline
50  &  157.08599750231211  &  0.5772156649015328 & $ 6.519\times 10^{-24} $  \\ \hline
60  &  188.50086358429127  &  0.5772156649015328 & $ 2.871\times 10^{-25} $  \\ \hline
70  &  219.91603253410894  &  0.5772156649015328 & $ 1.771\times 10^{-26} $  \\ \hline
80  &  251.33139082491842  &  0.5772156649015328 & $ 1.180\times 10^{-27} $  \\ \hline
90  &  282.74687536370536  &  0.5772156649015328 & $ 2.181\times 10^{-29} $   \\ \hline
100  &  314.16244828586940  &  0.5772156649015328 & $ 2.861\times 10^{-29} $  \\ \hline

 \end{tabular}
 \end{center}

\bigskip

\bigskip

{\bf Appendix:}  Below is a simple  PARI/GP  script  checking \eqref{EMRS1}  to arbitrary  accuracy
declared by   ${\setminus {\tt  p}}  $  precision,  below it is 2222.  The output gives  agreement  between lhs and rhs of
\eqref{EMRS1}  up to the number of digits  given by precision.  It takes a fraction of a second  to get results.
For really  big  precision  the allocated memory 300 MB maybe not sufficient.

\begin{verbatim}

allocatemem(300000000)

\p 2222

Soldner=solve(x=1.4, 1.5, real(eint1(-log(x))));
tmp=log(Soldner);
ss=suminf(n=1,  tmp^n/(n*n!));
write("EMRS.txt", Euler+log(tmp)+ss);

\end{verbatim}

In the above script we  used  the  fact that logarithmic integral is related  to the exponential integral ${\rm Ei}(x)$,
see e.g. \cite[chapt. 5]{abramowitz+stegun}:
\bee
{\rm li}(x)={\rm Ei}(\log x),
\eee
where
\bee
{\rm Ei}(x)=\equiv
\left\{
\begin{array}{ll}
\displaystyle{-p.v.\int_{-x}^\infty \frac{e^{-t}}{t} dt,~~{\rm for}~~x>0.}\\
\displaystyle{ -\int_{-x}^\infty \frac{e^{-t}}{t} dt,~~{\rm for}~~   x<0 }
\end{array}
\right.
\label{Ei}
\eee
and principal  value is needed to avoid singularity of the integrand at $t=0$.   The logarithmic integral is not implemented  in
Pari  while   exponential integral  is implemented as \verb"eint1(x)".    We have obtained as the result  the number
$4.27\times 10^{-2235}$.  To  check  \eqref{EMRS2}  change last  lines to 

\begin{verbatim}
ss=suminf(n=1, (-1)^n*tmp^n/(2^(n-1.0)*n!)*
         sum(k=0, floor((n-1)*0.5),  1.0/(2.0*k+1.0)));
write("EMRS.txt", Euler+log(tmp)-sqrt(Soldner)*ss);
\end{verbatim}

As the output this time we obtained  $2.7328\times 10^{-2233}$.

\bigskip

The equation \eqref{fifth}  can be checked in Pari  using  the following  commands:

\begin{verbatim}
allocatemem(500000000)
\p 2222

tmp=sumalt(n=1, (-1)^(n-1)/(2*n*(2*n)!));
oo=[1];
c_i=intnum(u=1, [oo, I], cos(u)/u);
print(Euler+c_i-tmp);
\end{verbatim}

The explanations are needed: PARI contains the numerical  routine \verb"sumalt" for summing infinite alternating
series  in which extremely efficient algorithm  of Cohen,  Villegas  and Zagier \cite{Zagier-Cohen} is implemented;
\verb+oo=[1]+ denotes in Pari infinity;   \verb+intnum( )+ is the function
for  numerical integration and flag \verb+k*I+ (\verb+I+$=\imath$, i.e. $\imath^2=-1$)  tells the procedure  that
the  integrand is an oscillating  function  of the type $\cos(kx)$,  here $k=1$.  After a few  minutes  we obtained
$1.42335\times 10^{-2235}$.    This  result  shows the power of Pari's procedures:  the value of the cosine integral
at  1  is  indeed    calculated  numerically  without  using the expansion \eqref{kosinus}  and the value of $\gamma$
to avoid a vicious circle (tautology).

\bigskip

\bigskip

{\bf Acknowledgment:} I thank Jonathan  Sondow and  Wadim Zudilin  for e--mail  exchange and important  remarks.  

\bigskip
\bigskip


\end{document}